\font\medbf=cmbx10 scaled \magstep1
\font\bigbf=cmbx10 scaled \magstep2
\magnification=1100
\pageno=0
\input psfig.sty
\tolerance 9000
\def\L{{\bf L}}
\def\ddd#1{\buildrel ...\over{#1}}

\def\v{\vskip 1em}
\def\vs{\vskip 2em}
\def\vsk{\vskip 4em}

\def\sqr#1#2{\vbox{\hrule height .#2pt
\hbox{\vrule width .#2pt height #1pt \kern #1pt
\vrule width .#2pt}\hrule height .#2pt }}
\def\square{\sqr74}
\def\endproof{\hphantom{MM}\hfill\llap{$\square$}\goodbreak}
\def\Ups{\Upsilon}

\def\a{\alpha}

\def\tv{\hbox{Tot.Var.}}

\def\Rc{{\cal R}}
\def\R{I\!\!R}

\def\i{\item}

\def\O{{\cal O}}
\def\S{{\cal S}}
\def\U{{\cal U}}

\def\C{{\cal C}}

\def\ve{\varepsilon}
\def\n{\noindent}
\null
\vsk\vsk
\centerline{\bigbf On the Boundary Control of}
\vs
\centerline{\bigbf Systems of Conservation Laws}
\vsk\vs
\centerline{\it Alberto Bressan and Giuseppe Maria Coclite}
\vs
\centerline{S.I.S.S.A., Via Beirut 4, Trieste 34014 Italy.}
\v
\centerline{e-mail: bressan@sissa.it, ~coclite@sissa.it}

\vskip 3cm

\n{\bf Abstract.} 
{The paper is concerned with  
the boundary controllability of entropy weak solutions
to hyperbolic systems of conservation laws. We prove a general result
on the asymptotic stabilization of a system near a constant state.
On the other hand, we give an example showing that 
exact controllability in finite time cannot be achieved, in general.}
\vsk\vsk\vsk
\vfill\eject

\n{\medbf 1 - Introduction}  
\v
Consider an $n\times n$ system of conservation laws on a bounded interval:
$$u_t+f(u)_x=0\qquad\qquad t\geq 0,~~x\in\,]a,b[\,.\eqno(1.1)$$
The system is assumed to be strictly hyperbolic, each
characteristic field being either linearly degenerate or genuinely nonlinear 
in the sense of Lax [8].  
We shall also assume that all characteristic speeds are
bounded away from zero.  More precisely, let $f:\Omega\mapsto\R^n$
be a smooth map, defined on an open set $\Omega\subseteq\R^n$.
For each $u\in\Omega$, call
$\lambda_1(u)<\cdots<\lambda_n(u)$ the eigenvalues of the Jacobian
matrix $Df(u)$.   We assume that there exists a minimum speed $c_0>0$ and 
an integer $p\in\{1,\ldots,n\}$ such that
$$\cases{\lambda_i(u)<0\qquad &if\quad $i\leq p$,\cr
\lambda_i(u)>0\qquad &if\quad $i> p$,\cr}\eqno(1.2)$$
$$
\big|\lambda_i(u)\big|\geq c_0>0\qquad\qquad u\in\Omega.\eqno(1.3)$$
By (1.2), for a solution defined on the strip $t\geq 0,~~x\in\,]a,b[\,$,
there will be $n-p$ characteristics entering at the boundary
point $x=a$, and $p$ characteristics entering at $x=b$.
The initial-boundary value problem is thus well posed
if we prescribe $n-p$ scalar conditions at $x=a$ and $p$
scalar conditions at $x=b$ [11].  
See also [1, 2] for the case of general entropy-weak 
solutions taking values in the space $BV$ of functions with bounded variation.

In the present paper we study the effect of boundary conditions 
on the solution of (1.1) from the point of view of control theory.
Namely, given an initial condition
$$u(0,x)=\phi(x)\qquad\qquad x\in\,]a,b[\eqno(1.4)$$
with small total variation, we regard the boundary data as
{\it control functions}, and study the family of configurations
$$\Rc(T)\doteq \big\{u(T,\cdot)\big\}\subset \L^1\big([a,b]\,;~\R^n\big)
\eqno(1.5)$$ 
which can be reached by the system
at a given time $T>0$.
\v
Beginning with the simplest case, consider a strictly hyperbolic system
with constant coefficients:
$$u_t+Au_x=0,\eqno(1.6)$$
where $A$ is a $n\times n$ constant matrix, with real distinct eigenvalues
$$\lambda_1<\cdots<\lambda_p<0<\lambda_{p+1}<\cdots<\lambda_n\,.$$
Call $$\tau\doteq \max_i {b-a\over|\lambda_i|}$$
the maximum time taken by waves to cross the interval $[a,b]$.
In this case, it is easy to see that the reachable set in (1.5)
is the entire space: $\Rc(T)=\L^1$
for all $T\geq \tau$.  In other words, the system is completely controllable
after time $\tau$.  Indeed, 
for any $T\geq\tau$ and initial 
and terminal data $\phi, \psi\in\L^1\big([a,b];\R^n\big)$,
one can always find a solution of (1.4), defined on the rectangle
$[0,T]\times [a,b]$ such that
$$u(0,x)=\phi(x),\qquad\qquad u(T,x)=\psi(x)\qquad\qquad x\in [a,b].$$
Such solution can be constructed as follows. 
Let $l_1,\ldots,l_n$ and $r_1,\ldots,r_n$ be dual bases of right and left
eigenvectors of $A$ so that $l_i\cdot r_j=\delta_{ij}$.
For $i=1,\ldots,n$, let $u_i(t,x)$ be a solution to the scalar
Cauchy problem
$$u_{i,t}+\lambda_i u_{i,x}=0,$$
$$ u_i(0,x)=\cases{ l_i\cdot \phi(x)
\qquad &if\qquad $x\in [a,b]$,\cr
l_i\cdot \psi(x+\lambda_iT)
\qquad &if\qquad $x\in [a-\lambda_iT,~b-\lambda_iT]$,\cr
~~0\qquad &otherwise.\cr}$$
Then the restriction of 
$$u(t,x)=\sum_i u_i(t,x) r_i$$
to the interval $[0,T]\times [a,b]$ satisfies (1.6)
and takes the required
initial and terminal values.  Of course, this corresponds to the solution 
of an initial-boundary value problem, determined by
the $n$ boundary conditions
$$\left\{\eqalign{l_i\cdot u(t,a)&=u_i(t,a)\qquad\qquad i=p+1,\ldots,n,\cr
l_i\cdot u(t,b)&=u_i(t,b)\qquad\qquad i=1,\ldots,p.\cr}\right.$$
This result on exact boundary controllability has been extended in
[9, 10] to the case of general quasilinear systems of the form
$$u_t+A(u)u_x=0.$$
In this case, the existence of a solution taking the
prescribed initial and terminal values is obtained for all
sufficiently
small data $\phi,\psi\in \C^1$.


Aim of the present paper is to study analogous controllability properties
within the context of 
entropy weak solutions $t\mapsto u(t,\cdot)\in BV$.
For the definitions and basic properties of weak solutions
we refer to [4].
For general nonlinear systems, it is clear that a 
complete controllability result within the space $BV$ cannot hold.
Indeed, already for a scalar conservation law, it was proved in [3]
that
the profiles $\psi\in BV$ which can be attained at a fixed time $T>0$
are only those which satisfy the Oleinik-type conditions
$$\psi'(x)\leq {f'\big(\psi(x)\big)\over (x-a)f'' \big(\psi(x)\big)}
\qquad\qquad \hbox{for a.e.}~~x\in [a,b].$$

For general $n\times n$ systems, a complete characterization
of the reachable set $\Rc(T)$ does not seem possible,
due to the complexity of repeated wave-front interactions.
 
Our first result is concerned with stabilization
near a constant state.  
Assuming that all characteristic 
speeds are bounded away from zero, we show that
the system can be asymptotically stabilized
to any state $u^*\in\Omega$, with quadratic rate of convergence.
\v
\noindent{\bf Theorem 1.} {\it Let $K$ be a compact, connected 
subset
of the open domain
$\Omega\subset\,\R^n$.  
Then there exist constants $C_0,\delta,\kappa>0$ such that the following 
holds.
For every constant state 
$u^*\in K$ and every initial data $u(0)=\phi:[a,b]\mapsto K$ with
$\tv\{\phi\}<\delta$,
there exists an entropy weak
solution $\> u = u (t, x) \>$ of (1.1) such that, for all $t>0$,}
$$\tv\big\{u(t)\big\}\leq C_0\,e^{-2^{\kappa t}}\,,\eqno(1.7)$$
$$\big\Vert u(t,x)-u^*\big\Vert_{L^\infty} \,\leq C_0\,e^{-2^{\kappa
t}}.\eqno(1.8)$$ 
\v
The proof will be given in Section 2.
An interesting question is whether the constant
state $u^*$ can be exactly reached, in a finite time $T$.
By the results in [9], this is indeed the case if the initial
data has small $\C^1$ norm.
On the contrary,
in the final part of this paper, we show that
exact controllability in finite time
cannot be attained in general, if the initial data is
only assumed to be small in $BV$.

Our counterexample is concerned with a class of strictly hyperbolic,
genuinely nonlinear $2\times 2$
systems of the form (1.1).  More precisely, we assume
\v
\i{(H)} The eigenvalues $\lambda_i(u)$ of the Jacobian matrix 
$A(u)=Df(u)$ satisfy
$$-\lambda^*<\lambda_1(u)<-\lambda_*~<0~<\lambda_*<\lambda_2(u)<\lambda^*\,.
\eqno(1.9)$$
Moreover, the right eigenvectors $r_1(u)$, $r_2(u)$ satisfy the inequalities
$$D\lambda_1\cdot r_1 >0,\qquad\qquad D\lambda_2\cdot r_2>0,\eqno(1.10)$$
$$r_1\wedge r_2<0,\qquad r_1\wedge(Dr_1\cdot r_1)<0,\qquad
r_2\wedge(Dr_2\cdot r_2)<0.\eqno(1.11)$$
\v
\n A partucular system which satisfies the above assumptions 
is the one studied by DiPerna [7]:
$$\left\{ \eqalign{\rho_t+(u\rho)_x&=0\,,\cr
u_t+\left( {u^2\over 2}+{K^2\over \gamma-1}\rho^{\gamma-1}\right)_x&=0\,,\cr}
\right.$$ 
with $1<\gamma<3$.
Here $\rho>0$ and $u$ denote the density and the velocity of a gas, 
respectively. 

The last two inequalities in (1.11) imply that the rarefaction curves
(i.e.~the integral curves of the vector fields $r_1,r_2$)
in the $(u_1,u_2)$ plane turn clockwise (fig.~1).
In such case, the interaction of two shocks of the same family
generates a shock in the other family. 

\midinsert
\vskip 10pt
\centerline{\hbox{\psfig{figure=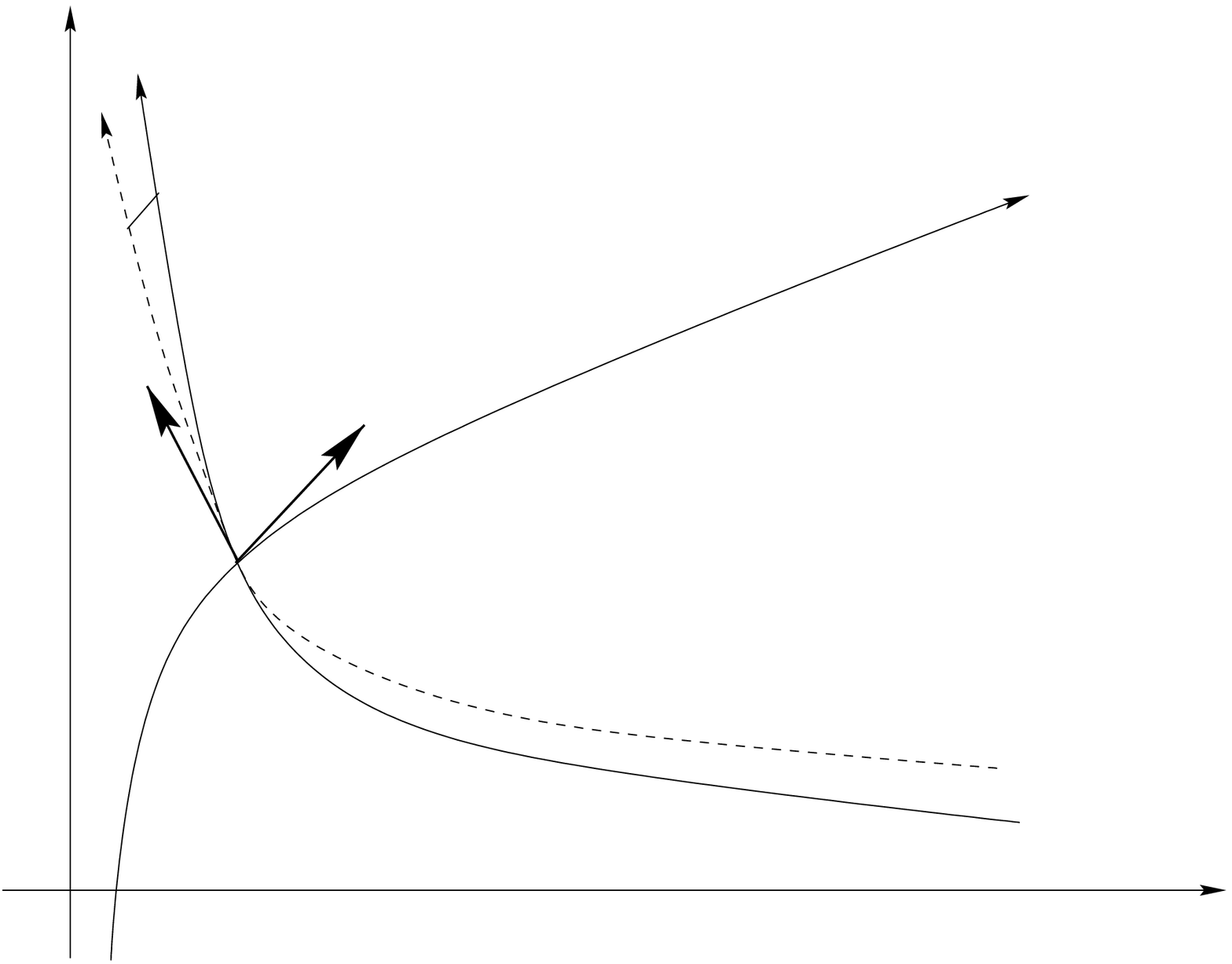,width=8cm}}}
\centerline{\hbox{figure 1}}
\vskip 10pt
\endinsert

\v
\n{\bf Theorem 2.} {\it Consider a $2\times 2$ system satisfying 
the assumption (H).  Then there exist initial data 
$\phi:[a,b]\mapsto\R^2$ having arbitrarily small total
bounded variation for which the following holds.
For every entropy weak solution $u$ of (1.1), (1.4),
with $\tv\big\{u(t,\cdot)\big\}$ remaining small for all $t$,
the set of 
shocks in $u(t,\cdot)$ is dense on $[a,b]$, for each $t>0$.   In particular,
$u(t,\cdot)$ cannot be constant.}
\v
As a preliminary,
in Section 3 we establish an
Oleinik-type estimate on the decay of positive waves.  This 
bound is of independent interest, and sharpens the results
in [5], for systems satisfying the additional conditions (H).

As a consequence, this implies that positive waves are
``weak'', and cannot completely cancel a shock within finite time.
The proof of Theorem 2 is then achieved by an induction argument. 
We show that, if the set of 1-shocks is dense on 
$[0,T]\times [a,b]$, then the set of points
$P_j=(t_j,x_j)$ where two 1-shocks interact and create a new 
2-shock is also dense on the same domain.
Therefore, new shocks are constantly generated, and the solution
can never be reduced to a constant. Details of the proof
will be given in Section 4.

As in [9], all of the above results refer to 
the case where total control on
the boundary values is available.  As a consequence, the problem is
reduced to proving the existence (or nonexistence) of an entropy weak
solution defined on the open strip  $t>0$, $x\in \,]a,b[\,$, satisfying the
required conditions.  This is a first step toward the analysis of
more general controllability problems, where the control
acts only on some of the
boundary conditions.  We thus leave open
the case where a subset of indices 
$I\subset\{1,\ldots,n\}$ is given, and one requires
$$l_i\cdot u(t,a)=\cases{\alpha_i(t)\qquad &if\quad $i\in I$,\cr
0\qquad &if\quad $i\notin I$,\cr}\qquad i=p+1,\ldots,n,$$
$$l_i\cdot u(t,b)=\cases{\alpha_i(t)\qquad &if\quad $i\in I$,\cr
0\qquad &if\quad $i\notin I$,\cr}\qquad i=1,\ldots,p,$$
for some control functions $\alpha_i$ acting only on the 
components $i\in I$.  
\vs
Throughout the following, we denote by $r_i(u)$, $l_i(u)$ the 
right and left $i$-eigenvectors of the Jacobian matrix $A(u)\doteq Df(u)$.
As in [4], we write
$\sigma \mapsto R_i (\sigma) (u_0)$ for the parametrized
$i$-rarefaction curve through the state $u_0$, so that
$${d \over {d\sigma}} R_i (\sigma)= r_i \big(R_i(\sigma)\big), 
\qquad R_i(0)= u_0.$$
The $i$-shock curve through $u_0$ is denoted by
$ \sigma \mapsto S_i (\sigma)(u_0)$.
It satisfies the Rankine-Hugoniot equations
$$f\big(S_i (\sigma)\big) - f(u_0) = 
\lambda_i(\sigma)\,\big(S_i(\sigma) - u_0 \big)$$
for some shock speed $\lambda_i$.
We recall (see [4], Chapter 5) that 
the general Riemann problem is solved in terms of the
composite curves
$$\Psi_i (u_0) (\sigma) =\cases{R_i (u_0) (\sigma),\quad
& ${\rm if}\>\> \sigma \geq 0,$\cr
S_i (u_0) (\sigma), & ${\rm if}\>\> \sigma <0. $\cr}\eqno(1.12)$$
\vfill\eject
\n{\medbf 2 - Proof of Theorem 1}  
\v
The proof relies on the two following two lemmas.
\v
\noindent{\bf Lemma 1.~~} {\it 
In the setting of Theorem 1, there exists a time $T>0$ such that
the following holds. For every pair of states
$\omega,\omega'\in K$
there exists an entropic
solution $\> u = u (t, x) \>$ of (1.1) such that} 
$$u ( 0, x ) \equiv \omega, \quad\quad\quad u ( T, x) \equiv 
\omega' 
\quad\quad\quad \hbox{for all~} x\in [a, b].\eqno(2.1)$$
\v
\noindent{\bf Proof.~} 
Consider the function
$$\Phi(\sigma_1,\ldots,\sigma_n\,;~v,\,v')\doteq
\Psi_n(\sigma_n)\circ\cdots\circ\Psi_{p+1}(\sigma_{p+1})(v')
-\Psi_p(\sigma_p)\circ\cdots\circ\Psi_1(\sigma_1)(v).\eqno(2.2)$$
Observe that, whenever $v=v'$, the $n\times n$ Jacobian matrix
$\partial \Phi/\partial \sigma_1\cdots\sigma_n$ computed at
$\sigma_1=\sigma_2=\cdots =\sigma_n=0$ has full rank.
Indeed, the columns of this matrix are given by the
linearly independent vectors $-r_1(v),\ldots,-r_p(v),\,r_{p+1}(v),\ldots,
r_n(v)$.
By the Implicit Function Theorem and a compactness argument we 
can find $\delta>0$ such that the following holds.
For every
$v,v' \in K$, with $|v-v'|\leq \delta$, there exist unique values
$\sigma_1,\ldots,\sigma_n$ such that
$$v''\doteq\Psi_n(\sigma_n)\circ\cdots\circ\Psi_{p+1}(\sigma_{p+1})(v')
=\Psi_p(\sigma_p)\circ\cdots\circ\Psi_1(\sigma_1)(v)\,.\eqno(2.3)$$

Defining the time 
$$\tau\doteq \max_{1\leq i\leq n}\sup_{u\in\Omega}{b-a\over \big|\lambda_i (u)
\big|}\,,\eqno(2.4)$$
we claim that there 
exists an entropy weak solution $u:[0,2\tau ]\times [a,b]\mapsto\Omega$
such that
$$u(0,x)\equiv v,\qquad\qquad u(2\tau,x)\equiv v'.\eqno(2.5)$$

\medskip
\midinsert
\vskip 10pt
\centerline{\hbox{\psfig{figure=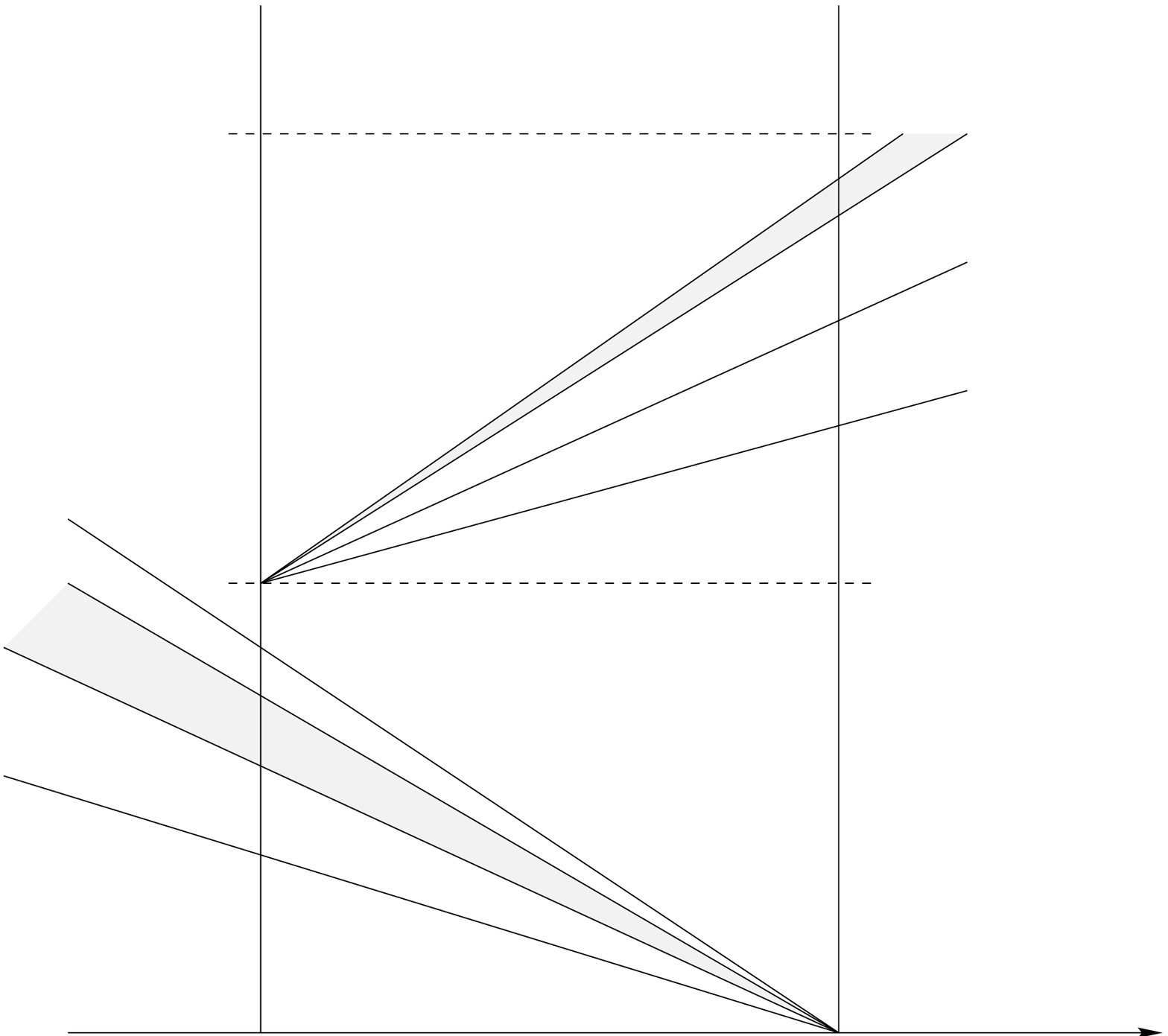,width=8cm}}}
\centerline{figure 2}
\vskip 10pt
\endinsert

The function $u$ is constructed as follows (fig.~2).
For $t\in [0,\tau]$ we let $u$ be the solution of the Riemann problem
$$u(0,x)=\cases{ v\qquad &if\quad $x<b$,\cr
v''\qquad &if\quad $x>b$.\cr}\eqno(2.6)$$
Moreover, for $t\in [\tau,~2\tau]$,
we define $u$ as the solution of the Riemann problem
$$u(\tau ,x)=\cases{ v'\qquad &if\quad $x<a$,\cr
v''\qquad &if\quad $x>a$.\cr}\eqno(2.7)$$
It is now clear that the restriction of
$u$ to the domain $[0,~2\tau]\times [a,b]$ satisfies the
conditions (2.5).  Indeed, by (2.3), on $[0,\tau]$
the solution $u$ contains only waves of families $\leq p$, 
originating at the point $(0,b)$. By (2.4) these waves cross the whole
interval $[a,b]$ and exit from the boundary point $a$ before time $\tau$.
Hence $u(\tau,x)\equiv v''$.
Similarly, still by (2.3), for $t\in [\tau,~2\tau]$ the function
$u$ 
contains only waves of families $\geq p+1$, 
originating at the point $(\tau,a)$. By (2.4) these waves cross the whole
interval $[a,b]$ and exit from the boundary point $b$ before time $2\tau$.
Hence $u(2\tau,x)\equiv v'$.

Next, given any two states $\omega,\omega'\in K$, 
by the connectedness assumption we can find a chain
of points $\omega_0=\omega,\omega_1,\ldots,\omega_N=\omega'$
in $K$ such that $|\omega_i-\omega_{i-1}|<\delta$ for every
$i=1,\ldots,N$.  Repeating the previous construction
in connection with each pair of states $(\omega_{i-1},\omega_i)$, 
we thus obtain
an entropy weak solution $u:[0,~2N\tau]\times [a,b] \mapsto \Omega $ 
that satisfies the conclusion of the lemma, with $T=2N\tau$.
\endproof
\v
In the following, we shall construct the 
desired solution $u=u(t,x)$ as limit of a sequence of front
tracking approximations.  Roughly speaking, an $\ve$-approximate
front tracking solution is a
piecewise constant function $u^\ve$, having jumps along a finite
set of straight lines in the $t$-$x$ plane say $x=x_\alpha(t)$, 
which approximately
satisfies the Rankine-Hugoniot equations:
$$\sum_\alpha \Big|f\big(u(t,x_\alpha+)\big)-f(u(t,x_\alpha-)\big)
-\dot x_\alpha \,\big(u(t,x_\alpha+)-u(t,x_\alpha-)\big)\Big|<\ve$$
for all $t>0$.   For details, see [4], p.125.
\v
\noindent{\bf Lemma 2.~~} {\it 
In the setting of Theorem 1, 
for every state $u^*\in \Omega$ there exist constants
$C,\delta_0>0$ for which the following holds.
For any $\ve>0$ and
every piecewise constant function
$\bar u:[a,b]\mapsto\Omega$ such that
$$\rho\doteq\sup_{x\in [a,b]}\big|\bar u(x)-u^*\big|\leq \delta_0,
\qquad\qquad
\delta\doteq \tv\{\bar u\}\leq\delta_0,\eqno(2.8)$$
there exists an $\ve$-approximate front tracking solution
$u=u(t,x)$ of (1.1), with $u(0,x)=\bar u(x)$, such that}
$$\sup_{x\in [a,b]}\big|u(3\tau,\,x)-u^*\big|\leq C\delta^2,
\qquad\qquad \tv\big\{u(3\tau)\big\}\leq  C\delta^2\,.\eqno(2.9)$$
\v
\noindent{\bf Proof. ~~} 
On the domain $(t,x)\in [0,\tau]\times [a,b]$, 
we construct $u$ as an $\ve$-approximate front
tracking solution in such a way that, whenever a front
hits one of the boundaries $x=a$ or $x=b$, no reflected front 
is ever created (fig.~3).
Since all fronts emerging from the initial data
$\bar u$ at time $t=0$ exit from $[a,b]$ within time $\tau$, 
it is clear that $u(\tau)$ can contain only fronts of
second or higher generation order. In other words, the only fronts 
that can be present in $u(\tau,\cdot)$ are the new ones,  
generated by interactions at times $t>0$ (the dotted lines in fig.~3). 
Therefore, using the interaction estimate (7.69) in [4] we obtain
$$\sup_{x\in [a,b]}\big|u(\tau,\,x)-u^*\big|=\O(1)\cdot 
(\rho+\delta)
\qquad\qquad \tv\big\{u(\tau)\big\}=\O(1)\cdot\delta^2\,.\eqno(2.10)$$

\midinsert
\vskip 10pt
\centerline{\hbox{\psfig{figure=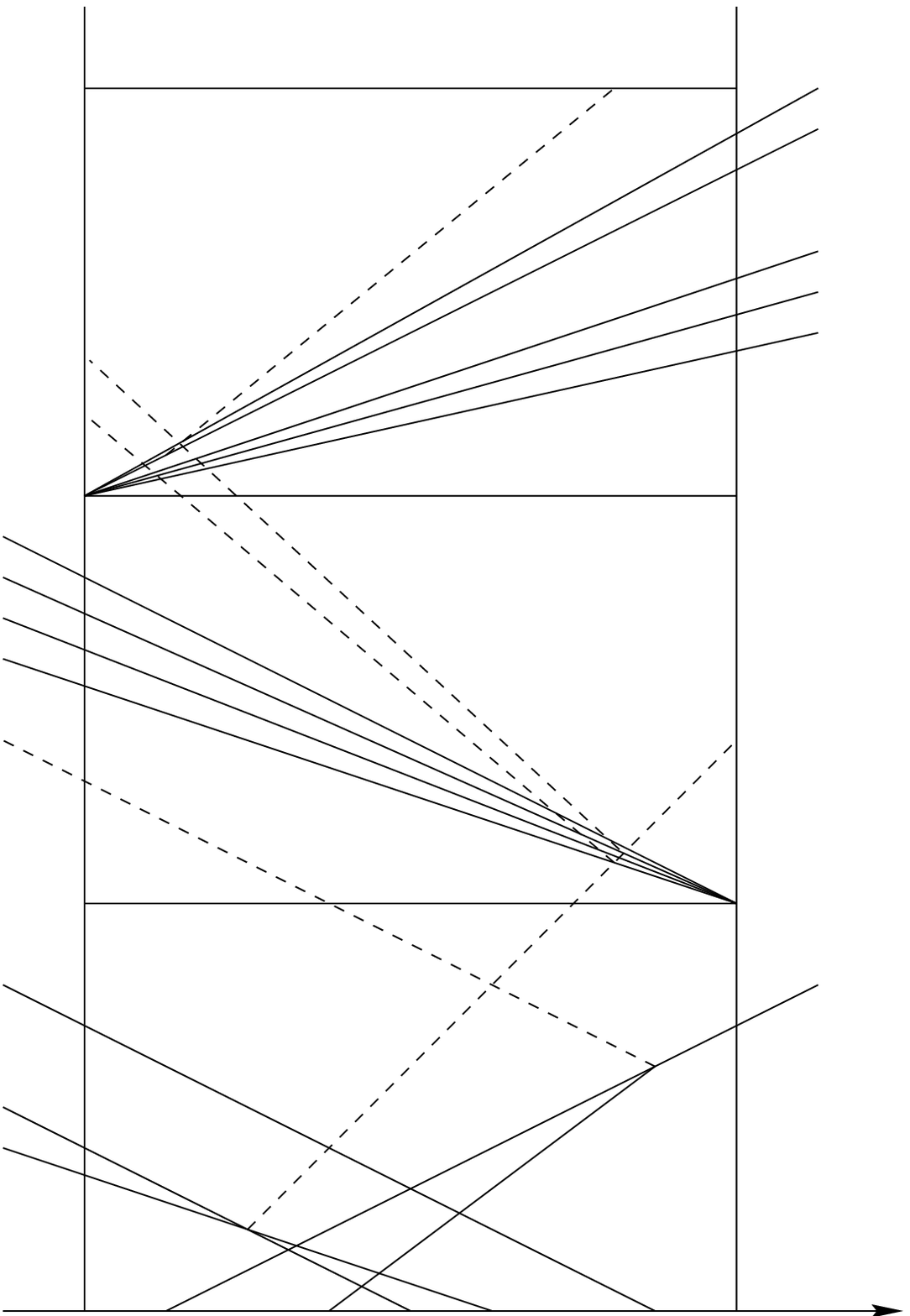,width=7cm}}}
\centerline{figure 3}
\vskip 10pt
\endinsert

\noindent We now apply a similar procedure as in the proof of Lemma 1, and
construct a solution on the interval $[\tau,~3\tau]$
in such a way that $u(3\tau)\approx u^*$. 
More precisely, to construct $u$ on the domain $[\tau,\,2\tau]\times [a,b]$,
consider the state $v''$ implicitly defined by
(2.2), with $v\doteq u(\tau, b-)$, $v'\doteq u^*$.
On a forward neighborhood of the point $(\tau,b)$ we let $u$ 
coincide with (a front-tracking approximation of) the solution
to the Riemann problem
$$u(\tau,x)=\cases{u(\tau,b-)\qquad &if\quad $x<b$,\cr
v''\qquad &if\quad $x>b$.\cr}$$
This procedure will introduce at the point $(\tau,b)$ 
a family of wave-fronts of families $i=1,\ldots,p$, whose total strength
is $\O(1)\cdot (\rho+\delta)$.
Because of (2.4), all 
these fronts will exit from the boundary $x=a$ within time $2\tau$.
Of course, they can interact with
the other fronts present in $u(\tau,\cdot)$. 
In any case, the total strength
of fronts in $u(2\tau,\cdot)$ is still estimated as
$$\tv\big\{u(2\tau)\big\}=\O(1)\cdot\delta^2\,.\eqno(2.11)$$
Next, to define $u$ for $t\in [2\tau,\,3\tau]$,
consider the state $v'''$ implicitly defined by
$$\left\{ \eqalign{ u(2\tau, a+)&=\Psi_n(\sigma_n)\circ\cdots\circ
\Psi_{p+1}(\sigma_{p+1})(v'''),\cr
u^*&=\Psi_p(\sigma_p)\circ\cdots\circ
\Psi_1(\sigma_1)(v''').\cr}\right.\eqno(2.12)$$
On a forward neighborhood of the point $(2\tau,a)$ we let $u$ 
coincide with (a front-tracking approximation of) the solution
to the Riemann problem
$$u(2\tau,x)=\cases{u(2\tau,\,a+)\qquad &if\quad $x>a$,\cr
v'''\qquad &if\quad $x<a$.\cr}$$
This procedure introduces at the point $(2\tau,\,a)$ 
a family of wave-fronts of families $i=p+1,\ldots,n$, whose total strength
is $\O(1)\cdot (\rho+\delta)$.
Because of (2.4), all 
these fronts will exit from the boundary $x=b$ within time $3\tau$.
Of course, they can interact with
the other fronts present in $u(2\tau,\cdot)$. 
In any case, the total strength
of fronts in $u(3\tau,\cdot)$ is still estimated as
$$\tv\big\{u(3\tau)\big\}=\O(1)\cdot\delta^2\,.\eqno(2.13)$$
Moreover, the difference between the values $u(3\tau,x)$ and $u^*$
will be of the same order of the total strength of waves in $u(\tau,\cdot)$,
so that the first inequality in (2.9) will also hold.
\endproof
\v
\noindent{\bf Proof of Theorem 1.~~} 
Using the same arguments as in the proof of
Lemma 1.1, for every $\ve>0$ we can construct an
$\ve$-approximate front tracking solution $\> u =u (t, x) \>$ on 
$[0, \,2N\tau]\times [a,b]$ such that
$$\sup_{x\in [a,b]}\big| u(2N\tau, x)-u^*\big|=\O(1)\cdot\delta,
\qquad\qquad
\tv\big\{u(2N\tau)\}=\O(1)\cdot\delta\,.\eqno(2.14)$$
Choosing $\delta>0$ sufficiently small, we can assume that, in (2.14),
$\O(1)\cdot\delta<\delta_0 < 1/C$, the constant in
Lemma 2.   Calling $T\doteq 2N\tau$, we  can now repeat the construction
described in Lemma 2 on each interval $\big[T+3k\tau,~T+3(k+1)\tau\big]$.
This yields
$$\sup_{x\in [a,b]}\big|u(T+3k\tau,\,x)-u^*\big|\leq \delta_k,
\qquad\qquad \tv\big\{u(T+3k\tau)\big\}\leq  \delta_k\,,\eqno(2.15)$$
where the constants $\delta_k$ satisfy the inductive relations
$$\delta_{k+1}\leq C\delta_k^2.\eqno(2.16)$$
Choosing a sequence of $\ve$-approximate front tracking solutions
$u_\ve$ satisfying (2.15)-(2.16) and taking the limit as
$\ve\to 0$, we obtain an entropy weak solution $u$
which still satisfies the same estimates.
The bounds (1.7)-(1.8) are now a consequence of (2.15)-(2.16), 
with a suitable choice of the constants
$C_0,\kappa$.
\endproof
\vsk
\n{\medbf 3 - Decay of positive waves}  
\v
Throughout the following, we consider
a $2\times 2$ system of conservation laws
$$u_t+f(u)_x=0\,,\eqno(3.1)$$
satisfying the assumptions (H).
Following 
[6], p.~128, we construct a set of Riemann coordinates
$\>(w_1,w_2)\>$.  One can then choose
the right eigenvectors of $Df(u)$ so that
$$r_i(u)={\partial u\over\partial w_i}\,,\qquad\quad
{\partial \lambda_i\over\partial w_i}=D\lambda_i\cdot r_i>0
\qquad\qquad i=1,2.\eqno(3.2)$$ 
It will be convenient to perform most of the analysis 
on a special class of solutions: 
piecewise Lipschitz functions
with finitely many shocks and no compression waves.
Due to the geometric structure of the system,
this set of functions turns out to be 
positively invariant for the flow
generated by the hyperbolic system.  
We first derive several a priori estimate concerning these solutions,
in particular on the strength and location of the shocks.
We then observe that any $BV$ solution can be obtained as limit
of a sequence of piecewise Lipschitz solutions in our special class.
Our
estimates can thus be extended to general $BV$ solutions.
\v
\noindent{\bf Definition 1.~~} {\it We call $\U$ the set of all
piecewise Lipschitz functions $u:\R
\mapsto\R^2$ with finitely many jumps, such
that:

\i{(i)} at every jump, the corresponding Riemann problem is solved
only in terms of shocks (no centered rarefactions);

\i{(ii)} no compression waves are present, i.e.: $w_{i,x}(x)\geq 0$ 
at almost every $x\in\R$, $i=1,2$.}
\v
The next lemma establishes the forward invariance of the set $\U$.
\v
\noindent{\bf Lemma 3.~~} {\it Consider the $2\times 2$ system of
conservation laws (3.1), satisfying the assumptions
(H).  
Let $\> u= u (t,x) \>$ be the solution to a Cauchy
problem, with small total variation, satisfying
$u ( 0, \cdot ) \in \U$.
Then}
$$u ( t, \cdot ) \in \U\qquad\quad \hbox{for all} ~~t\ge 0.\eqno(3.3)$$ 
\v
\noindent{\bf Proof.~} We have to show that, as time progresses, 
the total number of shocks does not increase and no compression wave
is ever formed.  This will be the case provided that
\v
\i{(i)} The interaction of two shocks of the same family
produces an outgoing shock of the other family.
\v
\i{(ii)} The interaction of a shock with an infinitesimal 
rarefaction wave of the same family produces
a rarefaction wave in the other family.
\v
Both of the above conditions can be easily checked by analysing the
relative positions of shocks and rarefaction curves.
We will do this for the first family, leaving the verification of the
other case to the reader.

Call $\sigma\mapsto R_1(\sigma)$ the rarefaction
curve through a state $u_0$,parametrized so that
$$\lambda_1\big(R_1(\sigma)\big)=\lambda_1(u_0)+\sigma\,.$$
It is well known that the shock curve through $u_0$ has a second order
tangency with this rarefaction curve. Hence there
exists a smooth function $c_1(\sigma)$ such that the point
$$S_1(\sigma)\doteq R_1(\sigma)+c_1(\sigma){\sigma^3\over 6}r_2(u_0)$$
lies on this shock curve, for all $\sigma$ in a neighborhood of zero.
{}From the Rankine-Hugoniot equations it now follows
$$\chi(\sigma)\doteq \Big(f\big(R_1(\sigma)+c_1(\sigma)(\sigma^3/ 6)\,r_2(u_0)
\big)
-f(u_0)\Big)\wedge \Big(R_1(\sigma)+c_1(\sigma)(\sigma^3/ 6)\,r_2(u_0)
-u_0\big)=0\,.
\eqno(3.4)$$
Differentiating the wedge product (3.4) four times at $\sigma=0$ and denoting
derivatives with upper dots, we obtain
$$\eqalign{{d^4\chi\over d\sigma^4}(0)&=
4\big[\lambda_1(u_0)\ddd R_1(0)+2\ddot R_1(0)+\lambda_2(u_0)\,c_1(0)r_2(u_0) 
\big]\wedge \dot R_1(0)\cr
&\qquad +6\big[\lambda_1(u_0)\ddot R_1(0)+\dot R_1(0)\big]\wedge\ddot R_1(0)
+4\lambda_1(u_0)\dot R_1(0)\wedge\big[\ddd R_1(0)+c(0)r_2(u_0)\big]\cr
&=4\big(\lambda_2(u_0)-\lambda_1(u_0)\big)c_1(0)\,r_2(u_0)\wedge r_1(u_0)
+2(Dr_1\cdot r_1)(u_0)\wedge r_1(u_0)\cr
&=0\,.\cr}$$
Hence
$$c_1(0) ={(Dr_1\cdot r_1)\wedge r_1\over 
2(\lambda_2-\lambda_1)(r_1\wedge r_2)}
<0\,.\eqno(3.5)$$
By (3.5), the relative position of 1-shock and 1-rarefaction curves
is as depicted in fig.~1.  By the geometry of wave curves, the properties
(i) and (ii) are now clear.   Figure 4a illustrates the interaction
of two 1-shocks, while fig.~4b shows the interaction between
a 1-shock and a 1-rarefaction.  
By $u_l,u_m,u_r$ we denote the left, middle and right states
before the interaction, while $u_m'$ is 
the middle state after the interaction. 
In the two cases, the solution of the Riemann problem
contains a 2-shock and a 2-rarefaction, respectively.
\endproof

\midinsert
\vskip 10pt
\centerline{\hbox{\psfig{figure=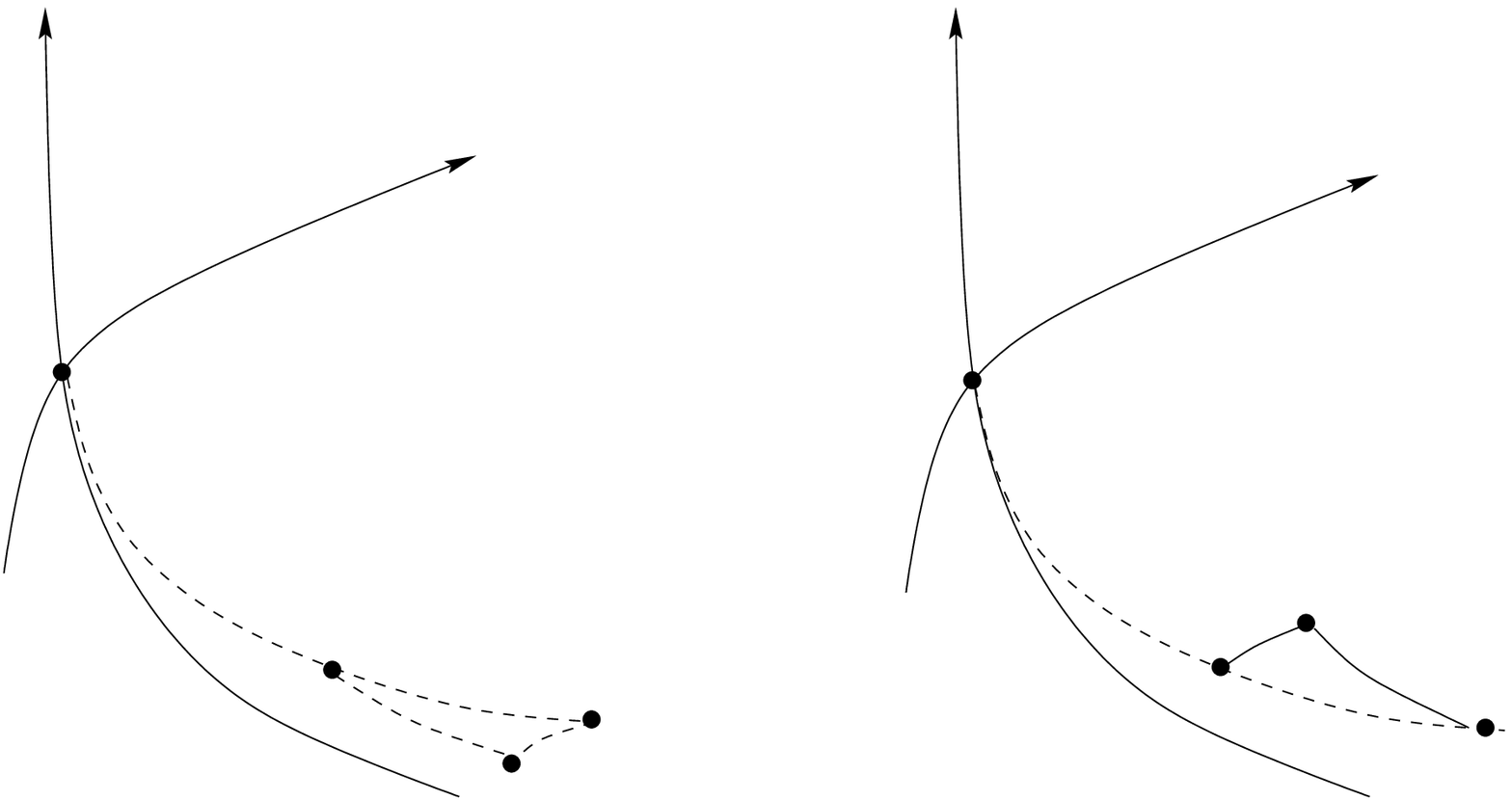,width=8cm}}}
\centerline{figure 4a~~~~~~~~~~~~~~~~~~~~~~~~~~~~~figure 4b}
\vskip 10pt
\endinsert

The next lemma shows the decay of positive waves for solutions
with small total variation, taking values inside $\U$.  
\v
\noindent{\bf Lemma 4.~~} {\it 
Let $\> u= u (t,x) \>$ be a solution of the Cauchy
problem for the
$2\times 2$ system (3.1) satisfying (H).  Assume that }
$$u ( t, \cdot ) \in \U \qquad\qquad t\ge 0. \eqno (3.6)$$ 
{\it Then there exist $\kappa, 
\delta >0 $ such that if $ \tv (u (t, \cdot ))<
\delta$ for all $t$, then its Riemann coordinates $(w_1,w_2)$
satisfy} 
$$ 0\leq w_{i,x}(t,x)\leq { \kappa \over t}, \qquad\qquad
t>0,~i=1,2.\eqno(3.7)$$
\v
\noindent{\bf Proof.~} We consider the case $ i =1$.  
Fix any point $(\bar t,~\bar x)$. Since
centered rarefaction waves are not present, there exists a unique
1-characteristic through this point, which we denote as
$\> t\mapsto x_1(t;~\bar t,~\bar x)$.  It is the solution of the Cauchy
problem
$$\dot x(t)=\lambda_1\big(u(t,x(t))\big),\qquad\qquad x(\bar t)=\bar
x.\eqno(3.8)$$ 
The evolution of $w_{1,x}$ along this characteristic is described by
$${d \over {dt}} w_{1,x} \big(t,x_1(t)\big) = w_{1,xt} + \lambda_1
w_{1,xx}=-( \lambda_1
w_{1,x} )_x + \lambda_1  w_{1,xx} = - {{\partial \lambda_1}
\over {\partial w_1}} w_{1,x}^2 - {{\partial \lambda_1}
\over {\partial \omega_2}} w_{1,x} w_{2,x}.$$
Since the
system is genuinely nonlinear there exists $k_1 >0$ such that
$\partial\lambda_1/\partial w_1\geq k_1 >0$, hence
$${d\over dt} w_{1,x} \big(t,x_1(t)\big)\leq -k_1 w^2_{1,x} 
+\O(1)\cdot w_{1,x}w_{2,x}. \eqno(3.9)$$
Moreover, at each time $t_\alpha$ where the characteristic crosses a 2-shock
of strength $|\sigma_\alpha|$
we have the estimate
$$w_{1,x}(t_\alpha+)\leq \big(1+\O(1)\cdot |\sigma_\alpha|\big)   
w_{1,x}(t_\alpha-).
\eqno(3.10)$$
Let $\> Q(t) \>$ be the total interaction potential at time $\> t \>$ (see
for example [4], p.~202) and  let 
$V_2(t) $ be the total amount of 2-waves approaching
our 1-wave located at $ x_1(t)$.  
Repeating the arguments in [4], p.139, we can find 
a constant
$C_0 >0$ such that the quantity 
$$\Upsilon(t)\doteq V_1(t)+C_0\,Q(t), \qquad\qquad t>0,$$
is non-increasing.
Moreover, for a.e.~$t$ one has
$$\dot \Ups (t)\leq - \big|\lambda_2-\lambda_1\big| 
|w_{2,x}|\big(t,x_1(t)\big)\,,$$
while at times $t_\alpha$ where $x_1$ crosses 
a 2-shock of strength $|\sigma_\alpha|$ there holds
$$\Ups(t_\alpha -)\leq \Ups(t_\alpha +)-|\sigma_\alpha|\,.$$
Call $\> W(t)\doteq w_{1,x} \big(t,x_1(t)\big)$. By the previous
estimates, from (3.9) and (3.10) it follows
form $$\dot W(t)\leq -k_1 W^2(t)- C\,\dot \Upsilon (t)W(t),\eqno(3.9)$$
$$W(t_\alpha+)-W(t_\alpha-)\leq 
C\big[ \Upsilon (t_\alpha+)-\Upsilon (t_\alpha-)\big]
W(t_\alpha-),\eqno(3.11)$$ 
for a suitable constant $C$.
We now observe that  $$y(t) \doteq {{ e^{-C \Upsilon(t)}} \over
\displaystyle{\int_0^t k_1 e^{-C \Upsilon(s)} ds}}$$
is a distributional solution of the equation
$$\dot y = -k_1 y^2 -C \,\dot \Upsilon (t) y\,,$$
with $y(t)\to \infty$ as $t\to 0+$.
A comparison argument now yields $W(t)\leq y(t)$.
Since $\> \Upsilon \>$ is positive and decreasing, we have 
$$W(t) \le \bar W (t) \le {1 \over k_1}{1 \over \displaystyle{\int_0^t e^{-C
\Upsilon(s)} ds}} \le {{e^{C\Upsilon(0)}} \over {k_1 t}},$$
for all $\> t>0$.   This establishes (3.7) for $i=1$, with 
$\kappa\doteq e^{C\Upsilon(0)}/k_1$.
The case $i=2$ is identical.
\endproof
\v
We conclude this section by proving a decay estimate
for positive waves, valid for general BV solutions
of the system (3.1).  For this purpose, we need to recall
some definitions introduced in [5]. See also p.~201 in [4].

Let $u:\R\mapsto \R^2$ have bounded variation. By possibly changing 
the values of $u$ at countably many points,
we can assume that $u$ is right continuous.
The distributional derivative $\mu\doteq D_x u$ is 
a vector measure,
which can be decomposed into a continuous and an atomic part:
$\mu=\mu_c+\mu_a$.  For $i=1,2$,
the scalar measures $\mu^i
=\mu^i_c+\mu^i_a$ are defined as follows.
The continuous part of $\mu^i$ is the Radon
measure $\mu^i_c$ such that
$$\int\phi~d\mu^i_c=\int \phi\,l_i(u)\cdot d\mu_c\eqno(3.12)$$
for every scalar continuous function $\phi$ with compact support.
The atomic part of $\mu^i$ is the measure $\mu^i_a$ concentrated on the
countable set $\{x_\a;~\a=1,2,\ldots\}$ where $u$ has a jump,
such that
$$\mu^i_a\big(\{x_\a\}\big)=\sigma_{\alpha,i}\doteq
E_i\big(u(x_\a-),~u(x_\a+)\big)\eqno(3.13)$$
is the size of the $i$-th wave in the solution of the
corresponding Riemann
problem with data $u(x_\a\pm)$. 
We regard $\mu^i$ as the {\it measure of $i$-waves} in the
solution $u$. It can be decomposed in a positive and a negative part,
so that 
$$\mu^i=\mu^{i+}-\mu^{i-},\qquad\qquad |\mu^i|=\mu^{i+} + \mu^{i-}.
\eqno(3.14)$$
The decay estimate in (3.7) can now be extended to general
BV solutions.  Indeed, we show that the density of positive 
$i$-waves decays as $\kappa/t$.
By meas$(J)$ we denote here the Lebesgue measure of a set $J$.
\v
\n{\bf Lemma 5.}~~{\it 
Let $u= u (t,x)$ be a solution of the Cauchy
problem for the
$2\times 2$ system (3.1) satisfying (H).
Then there exist $\kappa, 
\delta >0 $ such that if $ \tv (u (t, \cdot ))<
\delta$ for all $t$, then the measures $\mu^{1+}_t$, $\mu^{2+}_t$ 
of positive waves in $u(t,\cdot)$
satisfy
$$ \mu^{i+}_t(J)\leq { \kappa \over t}\,\hbox{meas}\,(J)\eqno(3.15)$$
for every Borel set $J\subset\R$ and every $t>0$, $i=1,2$.}
\v
\n{\bf Proof.}  For every $BV$ solution
$u$ of (3.1) we can construct a sequence of solutions $u_\nu$
with $u_\nu\to u$ as $\nu\to\infty$ and
such that $u_\nu(t,\cdot)\in\U$ for all $t$.
Calling $(w_1^\nu,\,w_2^\nu)$ the Riemann coordinates of $u_\nu$,
by Lemma 4 we have
$$ 0\leq w^\nu_{i,x}(t,x)\leq { \kappa \over t}, \qquad\qquad
t>0,~i=1,2,\qquad \nu\geq 1\,.\eqno(3.16)$$
For a fixed  $t>0$, observe that the map
$x\mapsto w^\nu_1(t,x)$ has upward jumps precisely at the
points $x_\alpha$ where $u(t,\cdot)$ has a 2-shock.
Define $\tilde\mu_\nu$ as the positive, 
purely atomic measure, concentrated on the
finitely many points $x_\alpha$ where $u(t,\cdot)$ has a 2-shock, such that
$$\tilde\mu_\nu
\big(\{x_\alpha\}\big)=w^\nu_1(t,\,x_\alpha +)-w^\nu_1(t,\,x_\alpha -)
\leq C\,|\sigma_\alpha|^3\eqno(3.17)$$
for some constant $C$.
By possibly taking a subsequence, we can assume the existence of a weak
limit $\tilde\mu_\nu\rightharpoonup \tilde\mu$.
Because of the estimate in (3.17), 
the measure $\tilde\mu$ is purely atomic, and
is concentrated on the set of points $x_\beta$ which are limits
as $\nu\to\infty$ 
of a sequence of points $x_\alpha^\nu$ where $u_\nu(t,\cdot)$ 
has a 2-shock of uniformly positive strength $|\sigma_\nu|\geq
\delta>0$.
Therefore, $\tilde\mu$ is concentrated on the set of points where
the limit solution $u(t,\cdot)$ has a 2-shock, and makes no contribution
to the positive part of $\mu^{1+}_t$.  We thus conclude that the
positive part of $\mu^{1+}_t$ is absolutely continuous w.r.t.~Lebesgue
measure, with density $\leq\kappa/t$.  An analogous argument holds for
$\mu^{2+}_t$.
\endproof
\v
\n{\bf Corollary 1.}  ~{\it 
Let $u= u (t,x)$ be a solution of the
$2\times 2$ system (1.1).  Let the assumptions (H) hold.
Fix $\ve>0$ and consider the subinterval $[a',b']\doteq [a+\ve,~b-\ve]$.
Assume that, at time $t=0$, the measures
$\mu^{1+}$, $\mu^{2+}$ 
of positive waves in $u(0,\cdot)$ on $[a,b]$
vanish identically.
Then, for every $t>0$ one has
$$ \mu^{i+}_t(J)\leq { \kappa \lambda^*\over \ve}\,
\hbox{meas}\,(J)\eqno(3.18)$$
for every Borel set $J\subset [a',b']$ and every $t>0$, $i=1,2$.}
\v
Indeed, recalling (1.9),
the values of $u(t,\cdot)$ restricted to the interval $[a',b']$
can be obtained by solving a Cauchy problem,
with initial data assigned on the whole
interval $[a,b]$ at time $t-\ve/\lambda^*$.
\vsk
\n{\medbf 4 - Proof of Theorem 2}
\v
\noindent{\bf Lemma 6.~~} {\it In the same setting as Lemma 4,
assume that there exists $\kappa'>0$ such that} 
$$ 0\leq w_{i,x}(t,x)\leq  \kappa'\qquad\qquad t\in [0,T],~~
i=1,2\,. \eqno(4.1)$$ 
{\it Let $\> t\mapsto x(t) \>$ be the location of a shock, with strength
$\big|\sigma(t)\big|$. There exists a constant $0< c< 1$ 
such that} 
$$\big|\sigma (t)\big|  \geq c \big| \sigma (s) \big|\,,\qquad
\qquad 0 \leq s <t\leq T\,.\eqno(4.2)$$
\v
\noindent{\bf Proof.~} 
To fix the ideas, let $u(t,\cdot)$ have a 1-shock located at $x(t)$,
with strength $\big|\sigma (t)\big|$. 
Outside points of interaction with other shocks, the strength
satisfies an inequality of the form
$${d \over {dt}} \big|\sigma (t) \big| \geq
-C \cdot \Big( w_{1,x} \big(t, x(t)+\big)+w_{1,x} \big(t, x(t)-\big)
w_{2,x} \big(t, x(t)+\big)+w_{2,x} \big(t, x(t)-\big)\Big)
\big|\sigma (t) \big|\,.\eqno(4.3)$$
At times where our 1-shock interacts with other 1-shocks, its strength
increases.  Moreover, at each time $t_\alpha$ where our 1-shock
interacts with a 2-shock, say of strength $|\sigma_\alpha|$,
one has
$$\big|\sigma(t_\alpha+)\big|\geq \big|\sigma(t_\alpha-)\big|\,
\big(1-C'|\sigma_\alpha|\big)\,.\eqno(4.4)
$$
for some constant $C'$.
Assuming that the total variation remains small, the 
total amount of 2-shocks which cross any given 1-shock
is uniformly small.  Hence, (4.3)-(4.4) together imply (4.2).
\endproof
\v
\noindent{\bf Lemma 7.~~} {\it 
Let $t\mapsto u(t,\cdot)\in\U$ be a solution of the Cauchy
problem for a genuinely
nonlinear $2\times 2$ system satisfying (1.11).  
Assume that there exists $\kappa'>0$ such that} 
$$ w_{i,x}(t,x)\leq  \kappa'\qquad\qquad t\in [0,T], ~~i=1,2\,. \eqno(4.5)$$ 
{\it 
Since no centered rarefactions are present, any two  
$i$-characteristics, say $x(t)<y(t)$,  can uniquely 
be traced backward up to time $ t=0$. 
There exists a constant $L >0
$ such that}  
$$y(t)-x(t)\leq L\, \big(y(s)-x(s)\big)\qquad\quad 0\leq s<t\leq T\,.
\eqno(4.6)$$
\v
\noindent{\bf Proof.~}  Consider the case $i=2$. 
By definition, the characteristics are solutions of
$$\dot x (t) = \lambda_2 \big( u ( t, x(t))\big), \qquad \dot y (t) =
\lambda_2 \big(u ( t, y(t)) \big).$$
Since the characteristic speed $\lambda_2$ decreases across 2-shocks,
we can write
$$\dot y(t)-\dot x(t)\leq C\,\int_{x(t)}^{y(t)}
\big|w_{1,x}(t,\xi)\big|+\big|w_{2,x}(t,\xi)\big|\,d\xi+
C\,\sum_{\alpha\in \S_1[x,y]}\big|\sigma_\alpha(t)\big|\,,\eqno(4.7)$$
where $\S_1[x,y]$ denotes the set of all 1-shocks located
inside the interval $\big[x(t),~y(t)\big]$.
Introduce the function
$$\phi(t,x)\doteq\cases{~~0\qquad &if\qquad $x\leq x(t)$,\cr
{x-x(t)\over y(t)-x(t)}\qquad &if $x(t)<x<y(t)$,\cr
~~1\qquad &if\qquad $x\geq y(t)$.\cr}$$
Moreover, define the functional
$$\Phi(t)\doteq \sum_{\alpha\in\S_1} \phi\big(t,x_\alpha(t)\big)\,
\big|\sigma_\alpha(t)\big|+ C_0\,Q(t)\,,$$
where the summation now refers to all 1-shocks in $u(t,\cdot)$ and
$Q$ is the usual interaction potential.
Observe that the map $t\mapsto \Phi(t)$ is non-increasing.
{}By (4.5) and (4.7) we can now write
$$\dot y(t)-\dot x(t)\leq C'\,\big(1-\dot\Phi(t)\big)\, 
\big( y(t)- x(t)\big)$$
for some constant $C'$.  This implies (4.6) with 
$L=\exp\big\{C'T+C'\Phi(0)\big\}$.
\endproof
\v
The next result is the key ingredient toward the proof of Theorem 2. 
It provides the density of the set of interaction points where
new shocks are generated.
\v
\noindent{\bf Lemma 8.~~} {\it Fix $\ve>0$ and define $a''=a+2\ve$,
$b''=b-2\ve$.  Consider a $2\times 2$ system of the form (1.1),
satisfying (H). 
Let $u$ be an entropy weak solution
defined on $ [0,\tau]\times [a,b]$, with 
$\tau\doteq \ve /4\lambda^*$. Let (3.18) hold for all $t\in [0,\tau]$, 
and assume that 
$u(0,\cdot)$ has
a dense set of 1-shocks on the interval $[a'',\,b'']$.
Then, for $0\leq t\leq \tau$, the solution 
$u(t,\cdot)$ has a set of 1-shocks
which is dense on $[a'',~b'-\lambda^*t]$ and a set
of 2-shocks
which is dense on $[a'',\,b'']$.}
\v
\noindent{\bf Proof.~} 
By the assumptions of the lemma, there exists
a sequence of piecewise  Lipschitz solutions $t\mapsto u_\nu(t) \in \U$
such that $u_\nu\to u$ in $\L^1$,
$$0\leq w_{i,x}^\nu(t,x)\leq {2\kappa\lambda^*\over \ve}
\qquad\qquad i=1,2,~~\nu\geq 1\,,$$
and moreover the following holds. 
For every $\rho>0$,
there exists $\delta>0$ such that each $u_\nu(0,\cdot)$ (with
$\nu$ large enough) contains at least one 1-shock of strength
$\big|\sigma_\nu(0)\big|\geq\delta$ 
on every subinterval $J\subset [a'',\,b'']$ having
length $\geq\rho$. 

To prove the first statement in Lemma 8, fix $t\in [0,\tau]$
and consider any non-trivial interval $[p,q]\subseteq [a'',~b''-t\lambda^*]$.
Call $s\mapsto p_\nu(s)$, $s\mapsto q_\nu(s)$ the backward
characteristics through these points, relative to the solution 
$u_\nu$.  We thus have
$$\left\{\eqalign{ \dot p_\nu(s)&=\lambda_1\big(u_\nu(s,\,p_\nu(s))\big),\cr
\dot q_\nu(s)&=\lambda_1\big(u_\nu(s,\,q_\nu(s))\big),\cr}
\right.\qquad\qquad\left\{\eqalign{p_\nu(t)&=p,\cr q_\nu(t)&=q.\cr}\right.$$
By Lemma 7, $q_\nu(0)-p_\nu(0)\geq \rho$ for some
$\rho>0$ independent of $\nu$. Hence, each solution $u_\nu$ contains
a shock of strength $\big|\sigma_\nu(s)\big|\geq\delta$ located 
inside the interval $\big[ p_\nu(0),~q_\nu(0)\big]$.
Lemma 5 now yields $\big|\sigma_\nu(t)\big|\geq c\delta$.
By possibly taking a subsequence, 
we conclude that the limit solution $u(t,\cdot)$
contains a 1-shock of positive strength at the 
point $x(t)=\lim x_\nu(t)\in [p,q]$.
\v
To prove the second statement, we will show that the set of points
where two 1-shocks in $u$ interact and produce a new 2-shock
is dense on the triangle 
$$\Delta\doteq \big\{ (t,x)\,;~~t\in [0,\tau], ~~a''<x<b''-\lambda^*t\big\}.$$
Indeed, let $t\in [0,\tau]$ and $p<q$ be as before. 
For each $\nu$ sufficiently large, let
$t\mapsto x_\nu(t)$ be the location of a 1-shock in $u_\nu$, 
with strength $\big|\sigma_\nu(t)\big|\geq\delta>0$.
Assume $x_\nu(\cdot)\to x(\cdot)$ as $\nu\to\infty$,
and $x_\nu(t)\in [p,q]$, so that $x(t)$ is the location of a 1-shock
of the limit solution $u$, say with strength $\big|\sigma(t)\big|>0$. 
\midinsert
\vskip 10pt
\centerline{\hbox{\psfig{figure=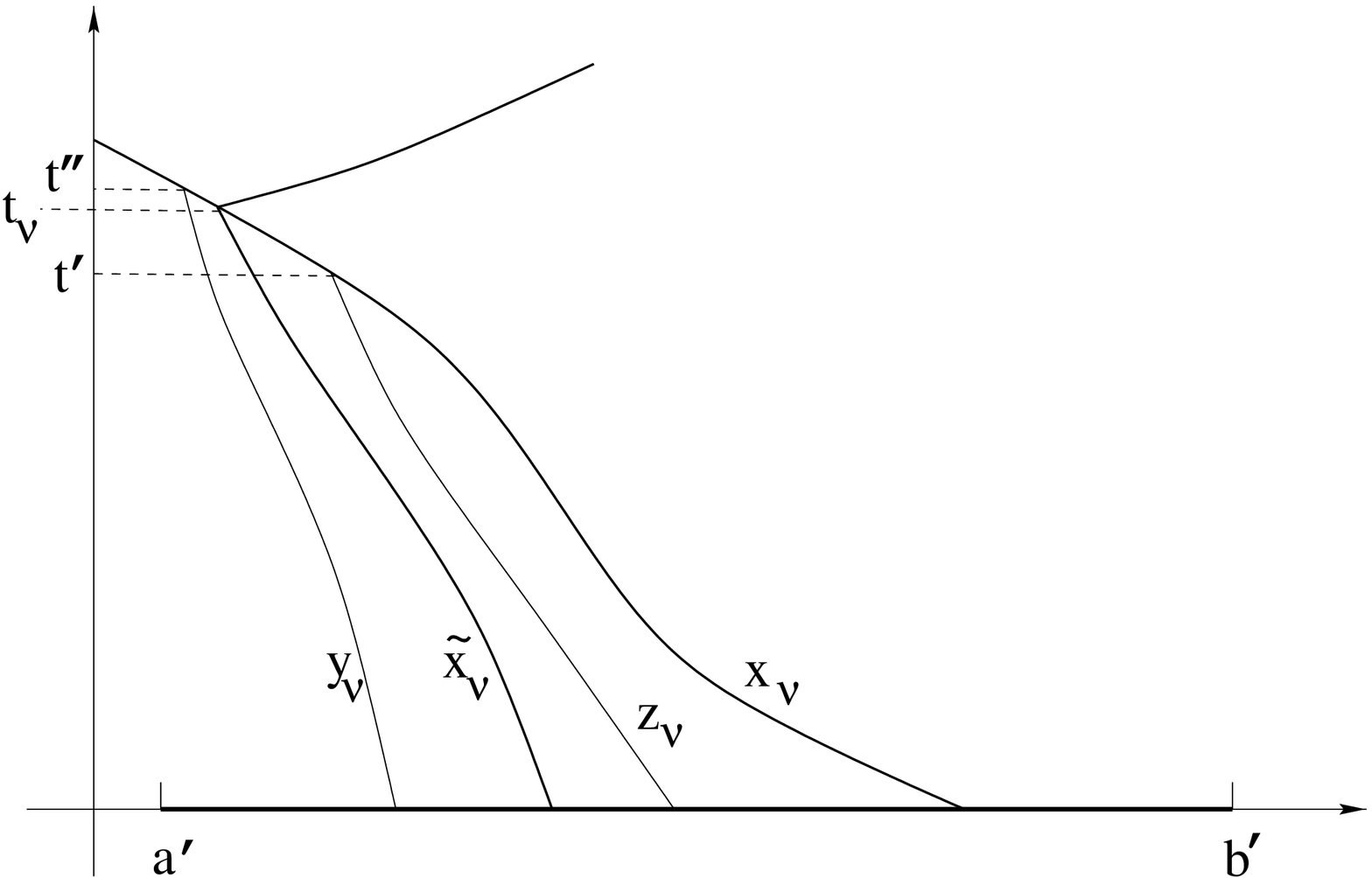,width=6cm}}}
\centerline{\hbox{figure 5}}
\vskip 10pt
\endinsert
We claim that the set of times $\hat t$ where some other 1-shock
$\sigma'$ impinges on $\sigma$ and generates a new 2-shock
is dense on $[0,t]$. To see this, fix $0<t'<t''<t$.
For each $\nu$ sufficiently large, 
consider the backward 1-characteristics $y_\nu,~z_\nu$
impinging from the left on the shock $x_\nu$
at times $t'',t'$ respectively (fig.~5).
These provide solutions to the
Cauchy problems $$\dot
y_\nu(t)=\lambda_1\big(u_\nu\big(t,y_\nu(t))\big),\qquad\qquad
y_\nu(t'')=x_\nu(t''),$$ $$\dot z_\nu(t)=
\lambda_1\big(u_\nu(t,z_\nu(t))\big),\qquad\qquad z_\nu(t')=x_\nu(t'),$$
respectively.
Observe that
$$z_\nu(0)-y_\nu(0)\geq\rho$$
for some $\rho>0$ independent of $\nu$. Indeed, the genuine
nonlinearity of the system implies
$$\lambda_1 \big( u_{\nu} (t, x_{\nu}(t)-) \big) - \dot x_{\nu} (t) 
\geq \kappa \Big|u_{\nu} \big(t, \,x_{\nu}(t)+) \big) - u_{\nu} \big(t, 
\,x_{\nu}(t)-) 
\Big| \geq \kappa \delta.$$
Therefore,
$$x_\nu(t')-y_\nu(t')\geq\rho'>0,$$
for some constant $\rho'>0$ independent of $\nu$.
By Lemma 6, the interval $\big[y_\nu(0),~z_\nu(0)\big]$ has uniformly
positive length.  Hence it contains a 1-shock of $u_\nu(0,\cdot)$
with uniformly positive 
strength $\big|\sigma_\nu(0)\big|\geq\delta>0$.
By Lemma 5, every $u_\nu$ 
has a 1-shock with strength $\big|\sigma_\nu(t)\big|\geq c\delta$
located along some curve
$t\mapsto \tilde x_\nu(t)$ with
$$y_\nu(t)<\tilde x_\nu(t)<z_\nu(t)\qquad t\in [0,t']\,.$$
Clearly, this second 1-shock impinges on the shock $x_\nu$ at
some time $t_\nu\in [t',t''],$ creating a new 2-shock with uniformly 
large strength.
Letting $\>\nu\to\infty\>$ we obtain the result.\endproof
\v
\noindent{\bf Proof of Theorem 2.~} 
Let $\delta_0>0$ be given.  We can then construct an initial condition
$u(0,\cdot)=\phi$, with $\tv\{\phi\}<\delta_0$,  
having a dense set of 1-shocks
on the interval $[a,b]$, and no other waves.
As a consequence, for any $\ve>0$ by Corollary 1 we 
have the estimate (3.18)
on the density of positive waves 
away from the boundary.

Fix $\tau=\ve/4\lambda^*$, and
consider again the subinterval $[a'',\,b'']=[a+2\ve,~b-2\ve]$.
We can apply Lemma 8 first on the time interval $[0,\tau]$, 
obtaining the density of 2-shocks on the region $[0,\tau]\times [a'',\,b'']$.
Then, by induction on $m$, the same argument is repeated on each
time interval 
$t\in \big[m\tau,~(m+1)\tau\big]$, proving the theorem.
\endproof
\vs
\n{\bf Acknowledgment.}  The second author 
warmly thanks professor Benedetto Piccoli for stimulating
conversations.
\vsk
\centerline{{\medbf References}}
\v
\i{[1]} D.~Amadori, Initial-boundary value problems for nonlinear systems of
conservation laws, {\it Nonlin. Diff. Equat. Appl.}  {\bf 4}
(1997), 1-42.
\v
\i{[2]} D.~Amadori and R.~M.~Colombo, Continuous dependence for
$2\times 2$ conservation laws with boundary, {\it J. Differential
Equations}  {\bf 138} (1997), 229-266.
\v
\i{[3]} F.~Ancona and A.~Marson, On the attainable set for scalar
nonlinear conservation laws with boundary control,
{\it SIAM J. Control Optim.} {\bf 36} (1998), 290-312.
\v
\i{[4]} A.~Bressan, {\it Hyperbolic Systems of Conservation Laws. 
The One Dimensional
Cauchy Problem}, Oxford University Press, 2000.
\v
\i{[5]} A.~Bressan and R.~M.~Colombo,
Decay of positive waves in nonlinear systems of conservation laws,
{\it Ann. Scuola Norm. Sup. Pisa} {\bf 
IV-26}
(1998), 133-160.
\v
\i{[6]} C.~M.~Dafermos, {\it 
Hyperbolic Conservation Laws in Continuum Physics}
Springer-Verlag, 2000. 
\v
\i{[7]} R.~DiPerna, Global solutions to a class of nonlinear hyperbolic 
systems of equations, {\it Comm. Pure Appl. Math.} {\bf 26} (1973),
1-28.
\v
\i{[8]} P.~Lax, Hyperbolic systems of conservation laws II, {\it
Comm. Pure Appl. Math.} {\bf 10} (1957), 537-566.
\v
\i{[9]} T.~Li and B.~Rao, Exact boundary controllability for quasilinear
hyperbolic systems, to appear.
\v
\i{[10]} T.~Li and J.~Yi, Semi-global $C^1$ solution
to the mixed initial-boundary value problem for quasilinear hyperbolic systems,
{\it Chin. Ann. of Math.} {\bf 21B} (2000), 165-186.
\v
\i{[11]} T.~Li and W.~Yu, 
{\it Boundary Value Problems for Quasilinear Hyperbolic Systems}, 
Duke University Mathematics Series V, 1985.

\bye